\magnification=1200

\input amstex

\documentstyle{amsppt}
%%%%%%%%%%%%%%%%%%%%%%%%%%%%%%%%%%%%%%%%%%%%%%%%%%%%%%%%%%%%%%

\hsize=165truemm

\vsize=227truemm

%%%%%%%%%%%%%%%%%%%%%%%%%%%%%%%%%%%%%%%%%%%%%%%%%%%%%%%%%%%%%%

\def\C{{{\Bbb C}}}

\def\p#1{{{\Bbb P}^{#1}_{\C}}}

\def\GL{\operatorname{GL}}

\def\PGL{\operatorname{PGL}}

\def\mapright#1{\mathbin{\smash{\mathop{\longrightarrow}
\limits^{#1}}}}

\def\mapdown#1{\Big\downarrow\rlap{$\vcenter{\hbox
{$\scriptstyle#1$}}$}}

%%%%%%%%%%%%%%%%%%%%%%%%%%%%%%%%%%%%%%%%%%%%%%%%%%%%%%%%%%%%%%

\topmatter

\title
On the rationality of moduli spaces of pointed curves 
\endtitle
\author
G. Casnati and C. Fontanari
\endauthor

\address
Gianfranco Casnati, Dipartimento di Matematica, Politecnico di Torino,
c.so Duca degli Abruzzi 24, 10129 Torino, Italy
\endaddress

\email
casnati\@calvino.polito.it
\endemail

\address
C. Fontanari, Dipartimento di Matematica, Politecnico di Torino,
c.so Duca degli Abruz\-zi 24, 10129 Torino, Italy
\endaddress

\email
claudio.fontanari\@polito.it
\endemail

\keywords
Pointed curve, Moduli space, Rationality
\endkeywords

\subjclassyear{2000}
\subjclass
14H10, 14H45, 14E08, 14L30
\endsubjclass

%\thanks
%This research was partially supported by MIUR and GNSAGA of INdAM 
%(Italy).
%\endthanks

\abstract 
Motivated by several recent results on the geometry of the moduli 
spaces $\overline{\Cal M}_{g,n}$ of stable curves of genus $g$ with $n$ 
marked points, we determine here their birational structure for small 
values of $g$ and $n$ by exploiting suitable plane models of the 
general curve. More precisely, we show that ${\Cal M}_{g,n}$ is 
rational for $g=2$ and $1\le n \le 12$, $g=3$ and 
$1\le n \le 14$, $g=4$ and $1 \le n \le 15$, $g=5$ 
and $1\le n\le 12$.
\endabstract

\endtopmatter

\document

\head
0. Introduction and notation
\endhead

The geometry of algebraic curves varying in families is a very old and fascinating subject, which naturally leads to consider the moduli space 
${\Cal M}_g$ of smooth projective curves of genus $g$ over the complex 
field $\C$. 

In order to investigate ${\Cal M}_g$ with the toolkit of algebraic 
geometry, the key ingredient is provided by its Deligne-Mumford 
compactification $\overline{\Cal M}_g$ defined in terms of stable 
curves. Indeed, $\overline{\Cal M}_g$ turns out to be a (mildly 
singular) projective variety and its boundary $\Delta=
\overline{\Cal M}_g \setminus {\Cal M}_g$ has pure 
codimension one. Moreover, $\Delta$ presents a natural inductive 
structure, its irreducible components being obtained from moduli 
spaces of strictly smaller genus via the operation of gluing curves 
at pairs of points. 

This remarkable fact suggests to consider more 
generally the moduli spaces ${\Cal M}_{g,n}$ of smooth $n$-pointed 
curves of genus $g$, parametrizing curves together with an ordered 
$n$-tuple of distinct points, and their analogous compactifications 
$\overline{\Cal M}_{g,n}$. Again, each irreducible component of 
the boundary $\overline{\Cal M}_{g,n} \setminus {\Cal M}_{g,n}$
is the image of either $\overline{\Cal M}_{g-1,n+2}$ or 
$\overline{\Cal M}_{g_1,n_1+1} \times \overline{\Cal M}_{g_2,n_2+1}$
with $g_1+g_2=g$ and $n_1+n_2=n$ via a suitable gluing morphism.
In order to exploit such an inductive pattern, the philosophy is to investigate these moduli spaces all together, 
even if one is mainly interested in $\overline{\Cal M}_{g}$.

Indeed, the inductive approach has recently reached significant results 
both from a topological and a geometrical point of view: we cannot resist 
to mention at least the computation of the first few rational cohomology 
groups of $\overline{\Cal M}_{g,n}$ by Arbarello and Cornalba in 
[A--C2] and the conjectural description of the ample cone of 
$\overline{\Cal M}_{g,n}$ in purely combinatorial terms by 
Gibney, Keel, and Morrison in [G--K--M]. 

Such a context also motivates the effort of extending to the pointed 
case the many techniques and results available for $\overline{\Cal M}_g$. For instance, the fact that its Kodaira dimension 
$\kappa(\overline{\Cal M}_g)$ is non-negative for every 
$g \ge 23$ has been established in a long series of papers (see [H--M], [Hr], [E--H], [Fa1]). Along these lines, a nice contribution [Lg]
by Logan determines for each $g$ with $4 \le g \le 23$ a small integer 
$h(g)$ such that the Kodaira dimension $\kappa(\overline{\Cal M}_{g,n})$ is non-negative for $n \ge h(g)$
(see also [Fa2] and [Fa3] for some recent improvements).

In the present paper instead we are going to prove the following

\proclaim{Main Theorem}
The moduli space
${\Cal M}_{g,n}$ is rational for $g=2$ and $1\le n \le 12$, $g=3$ and 
$1\le n \le 14$, $g=4$ and $1 \le n \le 15$, $g=5$ and $1\le n\le 12$.
\endproclaim

A few comments are in order here. 
The case $g=0$ is trivial: indeed, 
${\Cal M}_{0,n}$ is birational to $(\p1)^{n-3}$ for every $n \ge 3$, 
hence it is rational whenever it is defined. 
The case $g=1$ has been settled by P. Belorousski. In particular in [Be] it is proved via a nice geometric analysis of some plane cubic models of $n$--pointed elliptic curves that ${\Cal M}_{1,n}$ is rational for $1\le n\le 10$: for the sake of completeness we give a quick resum\'e of such a description in Section~1. We also point out that Belorousski's result is sharp: indeed, 
the fact that $H^{11,0}(\overline{\Cal M}_{1,11}, \C) \ne 0$ (see for 
instance [G--P], Section~4.1), implies that ${\Cal M}_{1,11}$ is not 
even unirational. More precisely, $\kappa(\overline{\Cal M}_{1,11})=0$ and $\kappa(\overline{\Cal M}_{1,n})=1$ for each $n\ge12$ (see [B--F], Theorem 3).

In higher genus the state of the art is as follows. 
It is well-known that ${\Cal M}_g$ is rational for $2\le g\le6$ (see respectively [Ig], [Ka4], [Ka2] and [SB1], [Ka3], [SB2]), and unirational for $7\le g\le14$ (see [A--S], [A--C1], [Se], [C--R], [Ve] and the references cited there).

The rationality of ${\Cal M}_{g,1}$ is also known for $g=2$ (see 
Section~6 of [Do]), $g=3$ (see Section~10 of [Do], where the 
rationality of ${\Cal M}_{3,1}$ is attributed to Shepherd-Barron 
without quoting any reference, and [Lo], Corollary~1.9, where an explicit proof is provided) and $g=4$ (see Theorem~B in [Ca]). 
Moreover, in Section~7 of [Lg] the author determines for 
$g=2, \ldots, 9, 11$ an integer $f(g)$ such that ${\Cal M}_{g,n}$ 
is unirational for $n \le f(g)$. 

The Main Theorem above shows that unirationality can be improved to 
rationality at least for genus up to $5$. 
We stress that for $g=2$ we even cover the first case in which 
unirationality was not previously known and for $g \ge 3$ 
our bounds for $n$ coincide exactly with the numbers $f(g)$ computed in [Lg], which turn out to be sharp at least in genus $g=4$. 
We also point out that our approach relies on classical constructions 
of suitable plane models of the general curve of genus $g$, 
obtained from its canonical embedding by using the marked points. 
Hence we deduce the existence of 
a birational equivalence ${\Cal M}_{g,n} \approx X / G$, where $G$ is a 
linear group and $X$ is a $G$-linearized projective bundle. Then we prove 
the rationality of ${\Cal M}_{g,n}$ by checking the rationality of the 
quotient $X / G$ via a case-by-case careful analysis exploiting several 
technical tools from rationality theory (for an overview of the standard 
methods see for instance [Do], [SB3], and Section~2.9 of [P--V]). Unfortunately, this natural strategy does not seem to apply for $g \ge 6$ and we leave the interested 
reader with the following

\proclaim{Main Question}
Is ${\Cal M}_{g,n}$ rational for all $n \le f(g)$, where $6 \le g \le 11$? 
\endproclaim

The authors wish to thank Gilberto Bini and Roberto Notari for several stimulating conversations on the subject of the present paper.

\subhead
Notation
\endsubhead
We work over the field $\C$ of complex numbers. We denote by $\GL_k$ the general linear group of $k\times k$ matrices with entries in $\C$ and by $\PGL_k$
the projective linear group, i.e. $\GL_k$ modulo the subgroup of scalar 
matrices.

Let ${\Bbb C}[x_1,\dots,x_k]$ be the ring of polynomials in the variables $x_1,\dots,x_k$ with coefficients in $\C$,  ${\Bbb C}[x_1,\dots,x_k]_d$ the
vector space of degree $d$ forms. 

If $g_1,\dots,g_h$ are elements of a certain group (resp. vector space) $G$ then $\langle
g_1,\dots,g_h\rangle$ denotes the subgroup (resp. subspace) of $G$ generated by $g_1,\dots,g_h$.

Projective $r$--space will be denoted by $\p r$: for each $i=1,\dots,r+1$ we denote by $E_i$ the point in $\p r$ all of whose coordinates but the $i^{th}$ one are zero: $E_{r+2}\in\p r$ is the unit point, i.e. the point whose coordinates are all $1$. If $V$ is a vector space then ${\Bbb P}(V)$ is the associated projective space.

If $X$ is a projective scheme, then a divisor $D$ on $X$ is always intended 
to be a Cartier divisor and $\vert D\vert$ denotes the complete linear system 
of divisors linearly equivalent to $D$.

A curve $C$ is a projective scheme of dimension $1$. If $D$ is a divisor on $C$ then $\varphi_{\vert D\vert}$ denotes the map induced on $C$ by $\vert D\vert$. If $C, D \subset \p2$ are distinct plane 
curves, then $C \cdot D$ denotes their intersection cycle. 

Let ${\Cal M}_{g,n}$ be the coarse moduli space of smooth, projective $n$--pointed curves (briefly $n$--pointed curves in what follows) of
genus $g$. A point of ${\Cal M}_{g,n}$ is then an ordered $(n+1)$--tuple $(C,p_1,\dots,p_n)$, where $C$ is a 
smooth and connected curve of genus $g$ and $p_1,\dots,p_n\in C$ are pairwise distinct point.
As usual, we set ${\Cal M}_g := {\Cal M}_{g,0}$.
A general curve of genus $g$ is a curve corresponding to the general 
point of the irreducible scheme ${\Cal M}_g$.

We denote isomorphisms by $\cong$ and birational equivalences by $\approx$.

For other definitions, results and notation we always refer to [Ht].

\head
1. The rationality of ${\Cal M}_{1,n}$ for $1\le n\le10$ (following Belorousski)
\endhead

In this section we recall some known results about the rationality of
${\Cal M}_{1,n}$ for $1\le n\le10$ due to P. Belorousski (see [Be]).
In order to follow Belorousski's construction, fix four points $A_1, 
A_2, A_3, A_4$ in $\p2$ in linearly general position and three lines 
$r_{1,2}$, $r_{2,4}$ and $r_{3,4}$ such that the line $r_{i,j}$ passes 
through the points $A_i$ and $A_j$.

\subhead
1.1. The case $n\le2$
\endsubhead
It is well-known that ${\Cal M}_{1,1} \cong \C$
via the $j$-invariant, so the first non-trivial case is $n=2$. 
In the space $\C[x_1,x_2,x_3]_3$ of cubic forms consider the subspace 
$V_2$ corresponding to cubics in $\p2$ passing through $A_1$, $A_2$, $A_3$, $A_4$ and tangent to $r_{1,2}$ at $A_1$, to $r_{2,4}$ at $A_2$, and to 
$r_{3,4}$ at $A_4$.
It is easy to see that ${ \Bbb P}(V_2) \cong \p2$ and that the general 
member of ${ \Bbb P}(V_2)$ is smooth. Hence there is a natural 
rational map $m_{1,2}: { \Bbb P}(V_2) \dashrightarrow {\Cal M}_{1,2}$ taking a cubic with its points $A_1$ and $A_2$ to the corresponding $2$-pointed curve of genus $1$. As checked in [Be], Lemma~1.1.4, $m_{1,2}$ turns out to be 
birational, so ${\Cal M}_{1,2}$ is rational.

\subhead
1.2. The case $3 \le n \le 4$
\endsubhead
Here the argument is very similar to the case $n=2$. Define ${ \Bbb P}(V_3) 
\subset \p9$ as the set of cubics in $\p2$ passing through $A_1$, 
$A_2$, $A_3$, $A_4$ and tangent to $r_{1,2}$ at $A_1$ and to $r_{2,4}$ 
at $A_2$,
and ${ \Bbb P}(V_4) \subset \p9$ as the set of cubics in $\p2$ passing through 
$A_1$, $A_2$, $A_3$, $A_4$ and tangent to $r_{1,2}$ at $A_1$. Exactly 
as in the previous subsection, there are natural isomorphisms  
${ \Bbb P}(V_3) \cong \p3$ and ${ \Bbb P}(V_4) \cong \p4$ and rational maps 
$m_{1,3}: { \Bbb P}(V_3) \dashrightarrow {\Cal M}_{1,3}$ and $m_{1,4}: 
{ \Bbb P}(V_4) \dashrightarrow {\Cal M}_{1,4}$ 
forgetting the plane embedding. Once again, in [Be] it is checked that 
both $m_{1,3}$ and $m_{1,4}$ are birational, hence the rationality of 
${\Cal M}_{1,3}$ and ${\Cal M}_{1,4}$ follows.

\subhead
1.3. The case $5 \le n \le 10$
\endsubhead
Let $V$ denote the vector space of cubic forms on $\p2$ vanishing at the 
points $A_1$, $A_2$, $A_3$, $A_4$. On each general cubic curve $C \in V$ 
the line $r_{1,2}$ cuts out a third point $A_5$ different from $A_1$ and 
$A_2$, thus the choice of four points uniquely determines a fifth point 
on the elliptic curve. Consider now the incidence variety
$$
X_n := \{\ (f, A_6, \ldots, A_n) \in{ \Bbb P}(V) \times (\p2)^{n-5} \ \vert\
f(A_i)=0,\ i=6,\dots,{n}\ \}.
$$
Since the dimension of $V$ is $6$, it follows that $X_n$ is non--empty 
for $n \le 10$. Moreover, Belorousski shows that $X_n$ is irreducible 
and its general point corresponds to a smooth curve (see Lemma~2.5.1 
below for a slightly weaker generalization).
Hence there is a canonical rational map $m_{1,n}: X_n \dashrightarrow 
{\Cal M}_{1,n}$ associating to a general point $(f, A_6, \ldots, A_n) \in X_n$ the genus 
one curve $\{ f=0 \}$ together with the $n$ distinct points $(A_1, \ldots, A_n)$. 
As proved by Belorousski, $m_{1,n}$ is indeed birational for all $n \ge 5$ 
(see [Be], Lemma~1.2.2) and $X_n$ turns out to be rational for $5 \le n \le 10$ (see [Be], Lemma~1.2.3). Hence ${\Cal M}_{1,n}$ is rational for all 
$n \le 10$.

\head
2. The rationality of ${\Cal M}_{2,n}$ for $1\le n\le12$
\endhead
Let $C$ be a general curve of genus $2$. Then $C$ is hyperelliptic and it carries exactly one $g^1_2$ which is the canonical series. If we choose four, possibly coinciding, points $p_1,p_2,p_3,p_4\in C$ then $D:=p_1+p_2+p_3+p_4$ defines a linear system $\vert D\vert$ which is a  $g_4^2$. By degree reasons $\vert D\vert$ has no base points, so we obtain a morphism $\varphi_{\vert D\vert}\colon C\to C_D\subseteq\p2$, whose image could be either an integral conic or a quartic divisor in $\p2$ with necessarily a double point by the genus formula for plane curves (see for instance [Ha], Example V.3.9.2). In the first case the $g_4^2$ should be composed with the unique canonical $g^1_2$, which is absurd if the points in the support of $D$ are general.

\remark{Remark 2.1}
If $\Psi\colon C\mapright\sim C'$ is an isomorphism of abstract curves  of genus $2$ such that $\Psi(p_i)=p_i'$, $i=1,\dots,4$ then $\Psi$ carries $\vert D\vert$ to $\vert D'\vert:=\vert p_1'+p_2'+p_3'+p_4'\vert$, thus it induces a projectivity $\psi\colon\p2\to\p2$ such that $\psi(C_D)=C'_{D'}:=\varphi_{\vert D'\vert}(C')$ making the following diagram
$$
\matrix
C&\mapright{\Psi}&C'\\
\\
\mapdown{\varphi_{\vert D\vert}}&&\mapdown{\varphi_{\vert D'\vert}}\\
\\
C_D&\mapright{\psi_{\vert C_D}}&C'_{D'}
\endmatrix\tag 2.1.1
$$
commutative. In particular $A_i:=\varphi_{\vert D\vert}(p_i)$ is mapped by $\psi_{\vert C_D}$ to $A_i':=\varphi_{\vert D'\vert}(p_i')$, $i=1,\dots,4$.

Conversely, each projectivity $\psi\colon \p2\to\p2$ mapping $C_D$ to $C'_{D'}$ and $A_i$ to $A_i'$, $i=1,\dots,4$,  is induced by an isomorphism $\Psi\colon C\mapright\sim C'$ such that $\Psi(p_i)=p_i'$, $i=1,\dots,4$. 
\endremark

\subhead
2.2. The case $n=1$
\endsubhead
If $(C,p_1)\in {\Cal M}_{2,1}$ is general then let $D:=4p_1$. The morphism 
$\varphi_{\vert D\vert}\colon C\to\p2$ maps $C$ to a quartic divisor $C_D\subseteq\p2$ with a $4$--flex and by the genus formula $C_D$ turns out to have also a double point. Up to projectivities we can assume that the singular point of $C_D$ is $E_1:=[1,0,0]$, the $4$--flex is $E_2:=[0,1,0]$ and the inflectional tangent line is $r:=\{\ x_1=0\ \}$. Thus the equation of $C_D$ is in 
$$
V:=x_1^2\C[x_2,x_3]_2\oplus x_1\C[x_2,x_3]_3\oplus \C x_3^4.
$$

Now the standard Bertini's theorem (see for instance [Ht], Theorem II.8.18 and Remark II.8.18.1) guarantees that the general element of $V$ corresponds
to an irreducible quartic divisor in $\p2$ with exactly one node at $E_1$.  
Hence we have a rational dominant map $m_{2,1}\colon V\dashrightarrow  {\Cal M}_{2,1}$ and we are going to identify its fibres. Remark~2.1 implies that $f,f'\in V$ are in the same fibre of $m_{2,1}$ if and only if there is a projectivity fixing $E_1$, $E_2$ and $r$ and sending $\{\ f=0\ \}$ to $\{\ f'=0\ \}$: thus it is represented by an element of the group
$$
G:=\left\{\ \pmatrix
\alpha_{1,1}&0&0\\
0&\alpha_{2,2}&\alpha_{2,3}\\
0&0&\alpha_{3,3}
\endpmatrix\in\GL_3\ \right\}.
$$
It follows that $ {\Cal M}_{2,1}\approx V/G$.

\subhead
2.3 The case $n=2$
\endsubhead
If $(C,p_1,p_2)\in {\Cal M}_{2,2}$ is general then we can repeat the argument for the $1$--pointed case with the divisor $D:=3p_1+p_2$. We can thus assume that the singular point of $C_D$ is $E_1$ and the flex is $E_2$ with inflectional tangent line $r:=\{\ x_1=0\ \}$ cutting $C_D$ also at $E_3:=[0,0,1]$. Thus the equation of $C_D$ is in 
$$
V:=x_1^2\C[x_2,x_3]_2\oplus x_1\C[x_2,x_3]_3\oplus \C x_2x_3^3.
$$

Once again, we have a birational equivalence $ {\Cal M}_{2,2}\approx V/G$, where
$$
G:=\left\{\ \pmatrix
\alpha_{1,1}&0&0\\
0&\alpha_{2,2}&0\\
0&0&\alpha_{3,3}
\endpmatrix\in\GL_3\ \right\}.
$$

\subhead
2.4 The case $n=3,4$
\endsubhead
If $(C,p_1,p_2,p_3)\in {\Cal M}_{2,3}$ (resp. $(C,p_1,p_2,p_3,p_4)\in {\Cal M}_{2,4}$) is general then we can repeat the argument for the $1$--pointed case with the divisor $D:=2p_1+p_2+p_3$ (resp. $D:=p_1+p_2+p_3+p_4$). We can thus assume that the singular point of $C_D$ is $E_1$ and that $p_1,p_2,p_3$ are respectively mapped to $E_2, E_3, E:=[0,1,1]$. Thus the equation of $C_D$ is in 
$$
V:=x_1^2\C[x_2,x_3]_2\oplus x_1\C[x_2,x_3]_3\oplus \C x_2(x_2-x_3)x_3^2
$$
(resp. in 
$$
V:=x_1^2\C[x_2,x_3]_2\oplus x_1\C[x_2,x_3]_3\oplus x_2x_3(x_2-x_3)\C[x_2,x_3]_1).
$$

Exactly as above, we obtain a birational equivalence $ {\Cal M}_{2,3}\approx V/G$ (resp. $ {\Cal M}_{2,4}\approx V/G$), where
$$
G:=\left\{\ \pmatrix
\alpha_{1,1}&0&0\\
0&\alpha_{2,2}&0\\
0&0&\alpha_{2,2}
\endpmatrix\in\GL_3\ \right\}.
$$

\subhead
2.5 The case $5\le n\le 12$
\endsubhead
If $(C,p_1,\dots,p_n)\in {\Cal M}_{2,n}$ is general then we can repeat the argument for the $4$--pointed case, so the equation of $C_D$ is again in
$$
V:=x_1^2\C[x_2,x_3]_2\oplus x_1\C[x_2,x_3]_3\oplus x_2x_3(x_2-x_3)\C[x_2,x_3]_1.
$$
Consider the incidence variety
$$
X:=\{\ (f,A_5,\dots,A_{n})\in {\Bbb P}(V)\times(\p2)^{n-4}\ \vert\ f(A_i)=0,\ i=5,\dots,{n}\ \}.
$$
Now we are going to define a rational dominant map $m_{2,n}\colon X\dashrightarrow{\Cal M}_{2,n}$ and to investigate its fibres. To this purpose we need the following standard result, which will be helpful also in the next sections.

\proclaim{Lemma 2.5.1}
Let $V\subseteq\C[x_1,x_2,x_3]_d$ be a subspace of dimension $v$, fix $m\le v-1$ and define
$$
X':=\{\ (f,A_1,\dots,A_m)\in {\Bbb P}(V)\times(\p2)^m\ \vert\ f(A_i)=0,\ i=1,\dots,m\ \}.
$$
Then there exists a maximal non--empty open subset ${\Cal U}\subseteq(\p2)^m$ such that $X:=X'\cap({\Bbb P}(V)\times{\Cal U})$ is a subbundle of ${\Bbb P}(V)\times{\Cal U}$ with fibre $\p{v-m-1}$ over ${\Cal U}$, in particular it is 
irreducible. Moreover the projection $X\to{\Bbb P}(V)$ is dominant.
\endproclaim
\demo{Proof}
Let $\{\ f_1, \ldots, f_v \ \}$ be a basis in $V$. Then each element $f\in V$ is of the form
$$
f(x_1,x_2,x_3)=\sum_{i=1}^v \alpha_i f_i(x_1, x_2, x_3).
$$
Let $x_{j,1},x_{j,2},x_{j,3}$ be coordinates on the $j^{th}$ copy of $\p2$ inside $(\p2)^m$. Then $X'\subseteq{\Bbb P}(V)\times(\p2)^m$ is defined by the equations
$$
\sum_{i=1}^v \alpha_i f_i(x_{j,1}, x_{j,2}, x_{j,3})=0,\qquad j=1,\dots,m.
$$
In particular, the fibres over $(\p2)^m$ are projective spaces of dimension at least $v-m-1$ linearly embedded inside ${\Bbb P}(V)$.

The dimension of the fibres over $(\p2)^m$ is an upper semicontinuous function and it is clear from the above equations that there are fibres of dimension exactly $v-m-1$, hence there is a non--empty open subset 
${\Cal U}\subseteq(\p2)^m$ over which the fibre is $\p{v-m-1}$. 
By restriction of the tautological section of ${\Bbb P}(V)\times{\Cal U}$, 
the subscheme $X:=X'\cap({\Bbb P}(V)\times{\Cal U})$ is indeed a subbundle,  
so it is irreducible of dimension $\dim (\p2)^m + \dim \p{v-m-1}=v+m-1$. 
Finally, since the fibres over a point of $X\to{\Bbb P}(V)$ are products of $m$ times the corresponding curve, the projection must be dominant.
\qed
\enddemo

We now come back to our case. Here the general element in ${\Bbb P}(V)$ is irreducible and smooth outside its unique double point $E_1$ (by Bertini's theorem again). It follows that the same is true for the general point in $X\approx \p{12-n}\times(\p2)^{n-4}$, thus we obtain a rational map $m_{2,n}\colon X\dashrightarrow{\Cal M}_{2,n}$. Let $G$ be the image inside $PGL_3$ of
$$
\left\{\ \pmatrix
\alpha_{1,1}&0&0\\
0&\alpha_{2,2}&0\\
0&0&\alpha_{2,2}
\endpmatrix\in\GL_3\ \right\}.
$$
Then the fibres of $m_{2,n}$ are the $G$--orbits of the points in $X$: in particular they have dimension $1$. Since $\dim(X)=n+4=\dim({\Cal M}_{2,n})+1$ we conclude that $m_{2,n}$ is dominant and, as in the previous cases, we have $ {\Cal M}_{2,n}\approx X/G$.

\medbreak
We finally conclude this section with 
\proclaim{Main Theorem for $g=2$} 
The moduli space $ {\Cal M}_{2,n}$ is rational for $1\le n\le12$.
\endproclaim
\demo{Proof}
On one hand, if $1\le n\le 4$ we have checked above that ${\Cal M}_{2,n}\approx V/G$ where $V$ is a linear representation of a solvable and connected algebraic group $G$. Thus the quotient on the right is rational by [Mi] or [Vi].

On the other hand, if $5\le n\le12$ we have checked above that $ {\Cal M}_{2,n}\approx X/G$ where $X$ is a $G$--linearized projective bundle over the rational base $(\p2)^{n-4}$ with typical fibre $\p{12-n}$.

The scheme $X$ is contained in the $G$--equivariant trivial projective bundle ${\Bbb P}(V)\times\Cal U$. It is not difficult to check the existence of a $G$--invariant unisecant $L\subseteq {\Bbb P}(V)\times\Cal U$  (i.e. a divisor intersecting the general fibre in a hyperplane). The scheme $L\cap X$ is then a $G$--invariant unisecant on $X$, hence $X$ is $G$--equivariantly birationally equivalent to a $G$--equivariant vector bundle $\C^{12-n}\times(\p2)^{n-4}$ for $5\le n\le12$.

Obviously, the action of $G$ on $(\p2)^{n-4}$ is almost free
(i.e., the stabilizer in $G$ of the general point of $(\p2)^{n-4}$ 
is trivial), thus ${\Cal M}_{2,n}\approx X/G\approx \C^{12-n}\times(\p2)^{n-4}/G$ is a vector bundle over the base $(\p2)^{n-4}/G$ with $(12-n)$--dimensional fibre by the results of Section~4 of [Do]. Moreover, $G$ is a torus and acts linearly on $(\p2)^{n-4}$, thus $(\p2)^{n-4}/G$ is rational (e.g. since $G$ is solvable and connected: see [Mi], [Vi]).
\qed
\enddemo

\head
3. The rationality of ${\Cal M}_{3,n}$ for $1\le n\le14$
\endhead
Let $C$ be a general curve of genus $3$. Since $C$ is not hyperelliptic, the canonical system $\vert K\vert$ is very ample, whence we obtain an embedding $\varphi_{\vert K\vert}\colon C\to C_K\subseteq\p2$. The image is a smooth quartic divisor.

\remark{Remark 3.1}
If $\Psi\colon C\mapright\sim C'$ is an isomorphism of abstract curves  of genus $3$ it carries  the canonical system of $C$ onto the one of $C'$. This is equivalent to saying that their canonical models $C_{K}$, $C'_{K'}$ are projectively equivalent. 
\endremark

\subhead
3.2. The case $n=1$
\endsubhead
Let $(C,p_1)\in {\Cal M}_{3,1}$. We can assume that $\varphi_K(p_1)=E_1$. Up to projectivities we can also assume that $C_K$ is tangent to the line $r:=\{\ x_3=0\ \}$ at the point $E_1$. Thus the equation of $C_K$ is in 
$$
V:=\C x_1^3x_3\oplus x_1^2\C[x_2,x_3]_2\oplus x_1\C[x_2,x_3]_3\oplus \C[x_2,x_3]_4.
$$

As in Section~2.2 there is a rational dominant map $m_{3,1}\colon V\dashrightarrow  {\Cal M}_{3,1}$ and again Remark~3.1 allows us to identify its fibres with the orbits of the natural action of the group
$$
G:=\left\{\ \pmatrix
\alpha_{1,1}&\alpha_{1,2}&\alpha_{1,3}\\
0&\alpha_{2,2}&\alpha_{2,3}\\
0&0&\alpha_{3,3}
\endpmatrix\in\GL_3\ \right\}
$$
on $V$. It follows that $ {\Cal M}_{3,1}\approx V/G$.

\subhead
3.3 The case $n=2,3$
\endsubhead
If $(C,p_1,p_2)\in {\Cal M}_{3,2}$ (resp. $(C,p_1,p_2,p_3)\in {\Cal M}_{3,3}$) we can assume that $\varphi_K(p_i)=E_i$, $i=1,2$ (resp. $i=1,2,3$), up to projectivities. Thus the equation of $C_K$ is in 
$$
V:=x_1^3\C[x_2,x_3]_1\oplus x_1^2\C[x_2,x_3]_2\oplus\langle x_1,x_3\rangle\C[x_2,x_3]_3
$$
(resp. in
$$
V:=x_1^3\C[x_2,x_3]_1\oplus \langle x_1^2,x_2x_3\rangle\C[x_2,x_3]_2\oplus x_1\C[x_2,x_3]_3).
$$

Once again, we have a birational equivalence $ {\Cal M}_{3,2}\approx V/G$ (resp. $ {\Cal M}_{3,3}\approx V/G$), where
$$
G:=\left\{\ \pmatrix
\alpha_{1,1}&0&\alpha_{1,3}\\
0&\alpha_{2,2}&\alpha_{2,3}\\
0&0&\alpha_{3,3}
\endpmatrix\in\GL_3\ \right\}
$$
(resp.
$$
G:=\left\{\ \pmatrix
\alpha_{1,1}&0&0\\
0&\alpha_{2,2}&0\\
0&0&\alpha_{3,3}
\endpmatrix\in\GL_3\ \right\}).
$$

\subhead
3.4 The case $4\le n\le14$
\endsubhead
If $(C,p_1,\dots,p_n)\in {\Cal M}_{3,n}$ we can assume that $\varphi_K(p_i)=E_i$, $i=1,2,3,4$, so that the equation of $C_D$ is in
$$
V:=\{\ f\in x_1^3\C[x_2,x_3]_1\oplus \langle x_1^2,x_2x_3\rangle\C[x_2,x_3]_2\oplus x_1\C[x_2,x_3]_3\ \vert\ f(1,1,1)=0\ \}.
$$
Let us consider the incidence variety
$$
X':=\{\ (f,A_5,\dots,A_{n})\in {\Bbb P}(V)\times(\p2)^{n-4}\ \vert\ f(A_i)=0,\ i=5,\dots,{n}\ \}
$$
and consider the irreducible open subscheme $X\subseteq X'$ and the open set ${\Cal U}\subseteq(\p2)^{n-4}$ defined in Lemma~2.5.1. Since the general element in ${\Bbb P}(V)$ is irreducible and smooth, the same is true for the general point in $X$. Since the unique automorphism of $\p2$ fixing the points $E_i$, $i=1,2,3,4$, is the identity we obtain that $ {\Cal M}_{3,n}\approx X\approx\p{14-n}\times(\p2)^{n-4}$ contained as a subbundle in ${\Bbb P}(V)\times(\p2)^{n-4}$.

\medbreak
Now the conclusion is analogous to the one of the previous section.
\proclaim{Main Theorem for $g=3$} 
The moduli space $ {\Cal M}_{3,n}$ is rational for $1\le n\le14$.
\qed
\endproclaim

\head
4. The rationality of ${\Cal M}_{4,n}$ for $1\le n\le15$
\endhead
Let $C$ be a general curve of genus $4$. Then the canonical system $\vert K\vert$ is very ample, whence we obtain an embedding $\varphi_{\vert K\vert}\colon C\to C_K\subseteq\p3$. The image is a sextic curve which is the complete intersection of a smooth quadric $Q$ with a cubic. Moreover two such sextic curves correspond to the same abstract curve if and only if they are projectively equivalent.

For each point $p_1\in C$ the divisor $D:=K-p_1$ defines a linear system $\vert D\vert$ which is a $g_5^2$. The map $\varphi_{\vert D\vert}\colon C\to C_D\subseteq \p2$ coincides with the projection of $C$ from $p_1$. 
Its image $C_D$ is a quintic divisor in $\p2$ passing through $A_1:=\varphi_{\vert D\vert}(p_1)$.

\proclaim{Proposition 4.1}
Under the above hypotheses, the curve $C_D$ is an integral quintic with two nodes $N_1,N_2$. The points $N_1,N_2$ and $A_1$ lie on a unique line $r$.

Conversely, each integral quintic in $\p2$ with exactly two nodes $N_1,N_2$ is obtained in the above way.
\endproclaim
\demo{Proof}
We start by proving the second part of the statement. We recall that the blow up of the plane at two points is isomorphic to the blow up of $Q$ at a point. Keeping this in mind, consider an integral quintic with two double points $N_1,N_2$ on the line $r$. Blowing up the plane at $N_1,N_2$ and blowing down the strict transform of $r$ we obtain a smooth curve of genus $4$ on $Q$ passing through a fixed point. The genus formula for curves on a smooth quadric then shows that such a curve is the canonical model $C_K$ of a suitable curve $C$ of genus $4$ with a fixed distinguished point $p_1$ and the restriction to $C_K$ of the blow up/down map coincides with $\varphi_{\vert K-p_1\vert}$.

Conversely, $C_D$ is obviously irreducible and it is also reduced by degree reasons. By the genus formula it must contain exactly two double points $N_1,N_2$. If $N_2$ is infinitely near to $N_1$, with a proper choice of the coordinates we can assume that $N_1:=[1,0,0]$ and that the tangent line at $N_1$ is $r:=\{\ x_3=0\ \}$. Then the equation of $C_D$ must lie in the space
$$
\C x_1^3x_3^2\oplus x_1^2x_3\C[x_2,x_3]_2\oplus x_1\C[x_2,x_3]_4\oplus\C[x_2,x_3]_5
$$
and a parameter computation shows that the family of plane quintics with such a singular point is irreducible of dimension $17$.

On the other hand, an analogous parameter computation shows that the family of plane quintics with two nodes is irreducible of dimension $18$. Taking into account the first part of the proof, it follows that for a general curve $C$ of genus $4$ and for a general point $p_1\in C$ the singularities of $C_D$ are two distinct nodes.

Finally, if $r$ is the line through $N_1,N_2$, then it is easy to check that $r\cdot C_D=2N_1+2N_2+A_1$.
\qed
\enddemo 

\remark{Remark 4.2}
If $\Psi\colon C\mapright\sim C'$ is an isomorphism of abstract curves  of genus $4$ mapping $p_1\in C$ to $p_1'\in C'$, then it carries  the canonical system $K$ of $C$ onto the one $K'$ of $C'$. Thus $\Psi$ maps $\vert D\vert:=\vert K-p_1\vert$ to $\vert D'\vert:=\vert K'-p_1'\vert$. 
Hence there exists a projectivity  $\psi\colon \p2\to\p2$ such that $\psi(C_{D})=C'_{D'}:=\varphi_{\vert D'\vert}(C')$ and satisfying a commutative diagram of the form (2.1.1).

Conversely, each projectivity $\psi\colon \p2\to\p2$ mapping $C_D$ to $C'_{D'}$ and $A_1$ to $A_1'$ is induced by an isomorphism $\Psi\colon C\mapright\sim C'$ such that $\Psi(p_1)=p_1'$. 
\endremark

\medbreak
The above Remark allows us to assume in what follows that the nodes are $N_i=E_i$, $i=1,2$, hence the line through the nodes is $r:=\{\ x_3=0\ \}$.
In particular, each projectivity $\psi\colon \p2\to\p2$ mapping $C_D$ to 
$C'_{D'}$ has to fix $\{\ E_1, E_2\ \}$.

\subhead
4.3. The case $n=1$
\endsubhead
If $(C,p_1)\in{\Cal M}_{4,1}$ and $D:=K-p_1$, then the equation of $C_D$ is in
$$
U:=x_1^3\C[x_2,x_3]_2\oplus \langle x_1^2,x_1x_3,x_3^2\rangle\C[x_2,x_3]_3.
$$

Moreover, the line $r:=\{x_3=0\}$ cuts out on $C_D$ the divisor $2E_1+2E_2+A_1$. By Bertini's theorem, there is a rational dominant map $m_{4,1}\colon U\dashrightarrow  {\Cal M}_{4,1}$ and Remark~4.2 allows us to identify its fibres with the orbits in $U$ of the natural action of the group
$H=H_0\langle\mu\rangle$, where
$$
H_0:=\left\{\ \pmatrix
\alpha_{1,1}&0&\alpha_{1,3}\\
0&\alpha_{2,2}&\alpha_{2,3}\\
0&0&\alpha_{3,3}
\endpmatrix\in\GL_3\ \right\},\qquad \mu:=\pmatrix
0&1&0\\
1&0&0\\
0&0&1
\endpmatrix.
$$
It follows that ${\Cal M}_{4,1}\approx U/H$.

Let $G:=T \langle \mu \rangle$, where $T$ is the torus of diagonal matrices in $\GL_3$ and
$$
V:=\langle x_1^3x_2^2,x_1^3x_3^2\rangle\oplus x_1x_3^2\C[x_2,x_3]_2\oplus \langle x_3^2,x_1^2\rangle\C[x_2,x_3]_3.
$$
We claim that $V$ is a $(H,G)$--section of $U$ in the sense of [Ka1]
(see also Section~2.8 of [P--V]). 

By definition, the first property to be checked is that $V$ is $G$-dense 
in $U$. Indeed, let $f\in U$ be of the form
$$
\aligned
f(x_1,x_2,x_3)&=a_1x_1^2x_2^3-2a_2x_1x_2^3x_3+b_1x_1^3x_2^2-2b_2x_1^3x_2x_3+\\
&+\text{monomials not in $\langle x_1^2x_2^3,x_1x_2^3x_3,x_1^3x_2^2,x_1^3x_2x_3\rangle$}.
\endaligned\tag4.3.1
$$
If we choose 
$$
h:=\pmatrix
1&0&a_2/a_1\\
0&1&b_2/b_1\\
0&0&1
\endpmatrix\in H
$$
then $h(f)\in V$ and we see that $V$ is $G$--dense in $U$. 

Next, we need to prove that if $f \in V$ is generic and $h \in H$ is such 
that $h(f) \in V$ then $h \in G$. Indeed, let $f\in V$ be as in formula (4.3.1): then $a_2=b_2=0$ and if $f$ is general enough we can suppose $a_1b_1\ne0$. If
$$
h:=\pmatrix
\alpha_{1,1}&0&\alpha_{1,3}\\
0&\alpha_{2,2}&\alpha_{2,3}\\
0&0&\alpha_{3,3}
\endpmatrix\in H_0
$$
then $h(f)$ contains the monomials $2a_1\alpha_{1,1}\alpha_{1,2}\alpha_{2,2}^3x_1x_2^3x_3+2b_1\alpha_{1,1}^3\alpha_{2,2}\alpha_{2,3}x_1^3x_2x_3$, hence it follows that $h(f)\in V$ if and only if $h\in T$. In a similar way one can argue for $h\in H_0\mu$, so the claim is checked.

Now from the properties of sections (see Proposition~1.2 of [Ka1]) we deduce ${\Cal M}_{4,1}\approx U/H\approx V/G$.

\subhead
4.4 The case $n=2$
\endsubhead
If $(C,p_1,p_2)\in{\Cal M}_{4,2}$, let again $D:=K-p_1$ and consider $C_D\subseteq \p2$.  We can assume that $\varphi_D(p_2)=E_3$. The equation of $C_D$ is then in 
$$
V:=\langle x_1^3,x_2x_3^2\rangle\C[x_2,x_3]_2\oplus \langle x_1^2,x_1x_3\rangle\C[x_2,x_3]_3.
$$

Once again, we have a birational equivalence $ {\Cal M}_{2,2}\approx V/G$, with $G=T\langle\mu\rangle$, $T$ being as above the torus of diagonal matrices in $\GL_3$.

\subhead
4.5 The case $3\le n\le 15$
\endsubhead
If $(C,p_1,\dots,p_n)\in {\Cal M}_{4,n}$, then we can assume that 
$\varphi_D(p_i)=E_{i+1}$, $i=2,3$, so the equation of $C_D$ is in
$$
V:=\{\ f\in \langle x_1^3,x_2x_3^2\rangle\C[x_2,x_3]_2\oplus \langle x_1^2,x_1x_3\rangle\C[x_2,x_3]_3\ \vert\ f(1,1,1)=0\ \}.
$$
Consider the irreducible scheme $X$ and the non--empty open subset ${\Cal U}\subseteq(\p2)^{n-3}$ defined in Lemma~2.5.1, so that $X\approx\p{15-n}\times {\Cal U}$. 
The general element in ${\Bbb P}(V)$ is irreducible and smooth outside its double points $E_1$ and $E_2$, hence the same is true for the general point in $X$. Since the only automorphisms of $\p2$ fixing the set $\{E_1, E_2\}$ and the points $E_3, E_4$ are the identity and  $\mu$, it follows that $ {\Cal M}_{4,n}\approx X/\langle\mu\rangle$.

\medbreak
The following theorem concludes this section.
\proclaim{Main Theorem for $g=4$} 
The moduli space $ {\Cal M}_{4,n}$ is rational for $1\le n\le15$.
\endproclaim
\demo{Proof}
Consider the case $n=1$ described in Section~4.3. We have thus a birational map ${\Cal M}_{4,1}\approx V/G$. 
Now let $V_1:=\langle x_2^2x_3^3,x_2x_3^4,x_1^2x_3^3,x_1x_3^4 \rangle  
\subseteq V$. Since $V_1$ is $G$--invariant, there exists a decomposition $V=V_1\oplus V_2$ for a suitable $G$--invariant subspace $V_2\subseteq V$, in particular the projection $v\colon V\to V_1$ gives a vector bundle structure to $V$.

L\"uroth's theorem guarantees that $V_1/G$ is rational. We now compute the kernel $K$ of the action of $G$ on $V_1$. Consider a general $c(x_1,x_2,x_3):=c_1x_2^2x_3^3+c_2x_2x_3^4+c_3x_1^2x_3^3+c_4x_1x_3^4\in V_1$. If
$$
g:=\pmatrix
\alpha_{1,1}&0&0\\
0&\alpha_{2,2}&0\\
0&0&\alpha_{3,3}
\endpmatrix\in T
$$
is in the stabilizer $G_c$ of $c$ inside $G$, then 
$\alpha_{2,2}^2\alpha_{3,3}^3=
\alpha_{2,2}\alpha_{3,3}^4=
\alpha_{1,1}^2\alpha_{3,3}^3=
\alpha_{1,1}\alpha_{3,3}^4=1$, thus $\alpha_{1,1}=\alpha_{2,2}=\alpha_{3,3}$ and $\alpha_{1,1}^5=1$. A similar computation shows that no elements of the form 
$$
g:=\pmatrix
0&\alpha_{1,1}&0\\
\alpha_{2,2}&0&0\\
0&0&\alpha_{3,3}
\endpmatrix\in T\mu
$$
can be in the stabilizer of $c$. We conclude that $K$ coincides with the normal subgroup of scalar fifth roots of the identity: in particular it coincides with the kernel of the action of $G$ on $V$. Thus  ${\Cal M}_{4,1}\approx V/G$ is rational by the results of Section~4 in [Do].

Now consider the case $n=2$ described in Section~4.4: once again, we have 
a birational map ${\Cal M}_{4,2}\approx V/G$ and letting $V_1:=\langle x_2^2x_3^3,x_2x_3^4,x_1^2x_3^3,x_1x_3^4\rangle\subseteq V$ we get a decomposition ${V}=V_1\oplus V_2$ into $G$--invariant subspaces. Thus we can prove the statement exactly as in the previous case.

Finally consider the case $n\ge3$ described in Section~4.5. In this case, $ {\Cal M}_{4,n}\approx X/\langle\mu\rangle$. The scheme $X\approx\p{15-n}\times{\Cal U}$ is a subbundle of ${\Bbb P}(V)\times{\Cal U}$ which obviously contains a $\langle\mu\rangle$--invariant unisecant $L$  (i.e. a divisor intersecting the general fibre in a hyperplane). Thus $L\cap X$ is still a $\langle\mu\rangle$--invariant unisecant on $X$. Hence there exists a $\langle\mu\rangle$--equivariant birational equivalence $X\approx\C^{15-n}\times(\p2)^{n-3}$ which is a $\langle\mu\rangle$--equivariant vector bundle. The group $\langle\mu\rangle$ acts on $(\p2)^{n-3}$ by exchanging the first two coordinates in each factor, so its action is trivially almost free. Moreover, this description gives us an obvious $\langle\mu\rangle$--equivariant birational equivalence $(\p2)^{n-3}\approx(\C^2)^{n-3}$. Once again, $\langle\mu\rangle$ exchanges the coordinates in each factor of $(\C^2)^{n-3}$, thus $(\C^2)^{n-3}$ is a linear representation of the abelian group $\langle\mu\rangle\cong{\Bbb Z}_2$ and the quotient is rational (the classical reference for this elementary fact 
is [Fi]). Now we conclude via the already quoted results in [Do].
\qed
\enddemo

\head
5. The rationality of ${\Cal M}_{5,n}$ for $1\le n\le12$
\endhead
Let $C$ be a general curve of genus $5$. Then the canonical system $\vert K\vert$ is very ample, whence we obtain an embedding $\varphi_{\vert K\vert}\colon C\to C_K\subseteq\p4$. The image is an octic curve which is the base locus of a net of quadrics.

For each pair of possibly coinciding general points $p_1,p_2\in C$ the divisor $D:=K-p_1-p_2$ defines a linear system $\vert D\vert$ which is a $g_6^2$. The map $\varphi_{\vert D\vert}\colon C\to C_D\subseteq \p2$ coincides with the projection of $C$ from the line joining $p_1$ and $p_2$ (if $p_1=p_2$ it coincides with the projection from the tangent line to $C$ at the point $p_1$). Its image $C_D$ is a sextic divisor in $\p2$ passing through the points $A_i:=\varphi_{\vert D\vert}(p_i)$.

\proclaim{Proposition 5.1}
Under the above hypotheses, the curve $C_D$ is an integral sextic with five nodes $N_1,\dots,N_5$. The points $N_1,\dots,N_5$ and $A_1,A_2$ lie on a unique integral conic $\widehat{C}$. If $A_1=A_2$ such a conic is tangent to $C$ at $A_1$.

Conversely, each integral sextic in $\p2$ with exactly five nodes $N_1,\dots,N_5$ lying on an integral conic $\widehat{C}$ is obtained in the above way.
\endproclaim
\demo{Proof}
We start by proving the second part of the statement.
Consider an integral sextic with five nodes $N_1,\dots,N_5$ in general position. Thus the blow up of the plane at the points $N_1,\dots,N_5$ is a quartic del Pezzo surface $\Sigma$. If we consider the anticanonical embedding of $\Sigma$ inside $\p4$, then the strict transform $\widetilde{C}$ of the curve is smooth and it is cut out on $\Sigma$ by a quadric. Since $\Sigma$ is the complete intersection of two quadrics, $\widetilde{C}$ is the smooth complete intersection of three quadrics in $\p4$ i.e. it is the canonical model $C_K$ of a suitable curve $C$ of genus $5$. The points $N_1,\dots,N_5$ lie on a unique integral conic $\widehat{C}$ cutting out only two other, possibly coinciding, points on our curve. The  strict transform of $\widehat{C}$ is a line on $\Sigma$ intersecting $C_K$ at two, possibly coinciding, points $p_1,p_2$. Then the restriction to $C_K$ of the blow up map coincides with $\varphi_{\vert K-p_1-p_2\vert}$.

Conversely, it is clear that $C_D$ is irreducible. In order to prove it is integral it then suffices to check it is reduced. If $C_D$ is a triple conic (resp. a double smooth cubic, a double singular cubic) then $C$ would be trigonal (resp. bielliptic, hyperelliptic), thus $C$ would not be general. 

Thus we can assume that $C_D$ is actually integral. Clearly $C_D$ is singular and by the genus formula it carries at most triple points as singularities. Since $C$, being general, is not trigonal, then $C_D$ carries only double points, which turn out to be exactly five (possibly 
infinitely near), again by the genus formula. We denote them by $N_1,\dots,N_5$. Since $C_D$ is integral such five nodes $N_1,\dots,N_5$ lie necessarily on a reduced conic $\widehat{C}$. If $\widehat{C}$ is integral then it is clearly unique. If $\widehat{C}$ is the union of two lines, one of them contains three of the five nodes $N_1,\dots,N_5$, the other one the remaining two points, hence it is unique also in this case. 

Sextics with five assigned double points on a reduced conic $\widehat{C}$ form a linear system of dimension $21$. If $\widehat{C}$ is reducible then the points $N_1,\dots, N_4$ can be chosen freely, the remaining one $N_5$ on the line joining $N_1$ and $N_2$. Thus the locus of sextics with a reducible $\widehat{C}$ is dominated by a projective bundle with fibre $\p{12}$ over the product $(\p2)^4\times\p1$, hence it is irreducible and its dimension is $21$. Analogously, other standard parameter computations show that each family of such plane integral
sextics with infinitely near double points has dimension at most $21$.

Similarly one can also check that the locus of plane integral sextics with five pairwise distinct nodes on an integral conic is irreducible of dimension $22$. Since we have proved above that each plane sextic with such a second configuration of singularities is the projection from a suitable secant line of its canonical model, it follows that for a general curve $C$ of genus $5$ and for general points $p_1,p_2\in C$ the singularities of $C_D$ are five nodes.

By construction, for each general line $\ell\subseteq \p2$ the divisor $\ell\cdot C_D+A_1+A_2$ is in the canonical system. Since the canonical system is cut out on $C_D$ by the adjoints of degree $3$ we have the linear equivalence on $C_D$
$$
\ell\cdot C_D+A_1+A_2\sim \ell\cdot C_D+\widehat{C}\cdot C_D-2\sum_{i=1}^5N_i.
$$
Thus $A_1+A_2\sim \widehat{C}\cdot C_D-2\sum_{i=1}^5N_i$ and since $C_D$ is not hyperelliptic we have
$\widehat{C}\cdot C_D=A_1+A_2+2\sum_{i=1}^5N_i$, i.e.
$A_1,A_2\in\widehat{C}$.
\qed
\enddemo

\remark{Remark 5.2}
As in the genus $4$ case, if $\Psi\colon C\mapright\sim C'$ is an isomorphism of abstract curves  of genus $5$ mapping $p_i\in C$ to $p_i'\in C'$, then it carries the canonical system $K$ of $C$ onto the one $K'$ of $C'$. Thus $\Psi$ maps $\vert D\vert:=\vert K-p_1-p_2\vert$ to $\vert D'\vert:=\vert K'-p_1'-p_2'\vert$. Hence there exists a projectivity  $\psi\colon \p2\to\p2$ sending $C_{D}$ to $C'_{D'}:=\varphi_{\vert D'\vert}(C')$ and inducing a commutative diagram of the form (2.1.1).

Conversely, each projectivity $\psi\colon \p2\to\p2$ mapping $C_D$ to $C'_{D'}$ and $A_i$ to $A_i'$ is induced by an isomorphism $\Psi\colon C\mapright\sim C'$ such that $\Psi(p_i)=p_i'$, $i=1,2$. 
\endremark

\medbreak
The above remark allows us to assume in what follows that $A_1=E_1$ (resp. $A_i=E_i$, $i=1,2$) and $\widehat{C}:=\{\ x_3^2-x_1x_2=0\ \}$. Notice that $\widehat{C}$ is the image of the Veronese embedding $v\colon\p1\to\p2$ given by $[s_1,s_2]\mapsto [s_1^2,s_2^2,s_1s_2]$. In particular, two divisors $C_{D},C'_{D'}\subseteq\p2$ turn out to be projectively equivalent if and only if they are mapped one onto the other by an automorphism of $\p2$ fixing $A_1$ (resp. $A_i$, $i=1,2$) and $\widehat{C}$. Every such automorphism is induced by an automorphism of $\widehat{C}\cong\p1$ fixing $v^{-1}(A_1)$ (resp. $v^{-1}(A_i)$, $i=1,2$), since it is well known restriction gives an isomorphism between the group of automorphisms of $\p2$ fixing $\widehat{C}$ and the group of automorphisms of $\widehat{C}\cong\p1$.

The five nodes correspond to an unordered $5$--tuple of pairwise distinct points $S_1,\dots,S_5$ in $\p1\setminus\{\ E_1\ \}$ (resp. $\p1\setminus\{\ E_1, E_2\ \}$). The set of such $5$--tuples can be naturally identified with the open subset ${\Bbb P}_5^0\subseteq {\Bbb P}_5:={\Bbb P}(\C[s_1,s_2]_5)$ corresponding to polynomials without multiple roots and not divisible by $s_2$ (resp. either $s_1$ or $s_2$).

\subhead
5.3. The case $n=1$
\endsubhead
If $(C,p_1)\in{\Cal M}_{5,1}$ and $D:=K-2p_1$, then $C_D$ is a sextic tangent to the line $\{\ x_2=0\ \}$ at the point $E_1$, hence it is defined by a polynomial in
$$
V:=\C x_1^5x_2\oplus x_1^4\C[x_2,x_3]_2\oplus x_1^3\C[x_2,x_3]_3\oplus x_1^2\C[x_2,x_3]_4\oplus x_1\C[x_2,x_3]_5\oplus \C[x_2,x_3]_6.
$$
If we consider the incidence variety
$$
Y:=\{\ (f,S_1,\dots,S_5)\in{\Bbb P}(V)\times{\Bbb P}_5^0\ \vert\ 
\text{$N_i:=v(S_i)$ is a node
of $\{\ f=0\ \}$, $i=1,\dots,5$}\ \}
$$
with projection $\pi_Y \colon Y\to{\Bbb P}_5^0$, then the above description shows that the equation of $C_D$ belongs necessarily to a fibre of $\pi_Y$. Conversely,

\proclaim{Lemma 5.3.1}
The scheme $Y$ is a subbundle of ${\Bbb P}(V)\times{\Bbb P}_5^0$ with typical fibre $\p{10}$ over ${\Bbb P}_5^0$, in particular it is irreducible.

Moreover the general point in $Y$ corresponds to an irreducible sextic curve tangent to the line $r$ at the point $E_1$, containing exactly five nodes on 
$\widehat{C}$ and no other singularities.
\endproclaim
\demo{Proof}
The first part of the above statement can be checked by imitating the proof of Lemma~2.5.1: in this case, using reducible curves it is easy to see that the dimension of the fibres of $\pi_Y$ over ${\Bbb P}_5^0$ is constant. The smoothness of the curve corresponding to the general point of $Y$ can be easily proved by applying Bertini's theorem to suitable reducible curves singular at the points $N_1,\dots,N_5$.
\qed
\enddemo

Thus there is a rational dominant map $m_{5,1}\colon Y\dashrightarrow  {\Cal M}_{5,1}$: Remark~5.2 allows us to identify its fibres with the orbits in $Y$ of the natural action of the stabilizer in $\GL_2$ of the point $[1,0]=v^{-1}(A_1)$, which is the image $H$ inside $PGL_2$ of
$$
\left\{\ \pmatrix
\alpha_{1,1}&\alpha_{1,2}\\
0&\alpha_{2,2}
\endpmatrix\in\GL_2\ \right\}.
$$
It follows that ${\Cal M}_{5,1}\approx Y/H$.

Let $\overline{\Bbb P}_5^0:={\Bbb P}(\langle s_1^5,s_1^3s_2^2,s_1^2s_2^3,s_1s_2^4,s_2^5\rangle)\cap{\Bbb P}_5^0$, and consider the incidence variety $X:=\pi_Y^{-1}( \overline{\Bbb P}_5^0)$ with projection $\pi:=(\pi_Y)_{\vert X}\colon X\to \overline{\Bbb P}_5^0$. 
The variety $X$ is generically a projective bundle with typical fibre $\p{10}$ over the rational base $B=\overline{\Bbb P}_5^0$.

Let $T\cong\C^*\subseteq\PGL_2$ be the image of the torus of diagonal matrices. We claim that $X$ is a $(H,T)$--section of $Y$. Since $X=\pi^{-1}_Y(\overline{\Bbb P}_5^0)$ is irreducible and it dominates $\overline{\Bbb P}_5^0$, 
by Proposition~1.2 of [Ka1] it is sufficient to check that $\overline{\Bbb P}_5^0$ is a $(H,T)$--section of ${\Bbb P}_5^0$: this last assertion is an easy exercise using Tschirnhaus transformations.

We conclude that ${\Cal M}_{5,1}\approx Y/H\approx X/T$.

\subhead
5.4 The case $n=2$
\endsubhead
If $(C,p_1,p_2)\in{\Cal M}_{5,2}$, set $D:=K-p_1-p_2$ and consider $C_D\subseteq \p2$. $C_D$ is a sextic through the points $E_1, E_2$, hence it is defined by a polynomial in
$$
V:=x_1^5\C[x_2,x_3]_1\oplus x_1^4\C[x_2,x_3]_2\oplus x_1^3\C[x_2,x_3]_3\oplus x_1^2\C[x_2,x_3]_4\oplus\langle x_1,x_3\rangle\C[x_2,x_3]_5.
$$
Once again, define the incidence variety 
$$
X:=\{\ (f,S_1,\dots,S_5)\in{\Bbb P}(V)\times{\Bbb P}_5^0\ \vert\ 
\text{$N_i:=v(S_i)$ is a node on $\{\ f=0\ \}$, $i=1,\dots,5$}\ \}.
$$
As in the previous case, the following result holds.

\proclaim{Lemma 5.4.1}
The scheme $X$ is a subbundle of ${\Bbb P}(V)\times{\Bbb P}_5^0$ with typical fibre $\p{10}$ over ${\Bbb P}_5^0$, in particular it is irreducible.

Moreover the general point in $Y$ corresponds to an irreducible sextic curve through the points $E_1,E_2$, carrying exactly five nodes on $\widehat{C}$ and no other singularities.
\qed
\endproclaim

As in Section~5.3, ${\Cal M}_{5,2}\approx X/T$.

\subhead
5.5 The case $3\le n\le 12$
\endsubhead
If $(C,p_1,\dots,p_n)\in {\Cal M}_{5,n}$, set once again $D:=K-p_1-p_2$ and consider $C_D\subseteq \p2$. Let
$$
V:=x_1^5\C[x_2,x_3]_1\oplus x_1^4\C[x_2,x_3]_2\oplus x_1^3\C[x_2,x_3]_3\oplus x_1^2\C[x_2,x_3]_4\oplus\langle x_1,x_2\rangle\C[x_2,x_3]_5
$$
and define the incidence variety
$$
\align
X':=\{\ &(f,S_0,\dots,S_4,A_3,\dots,A_n)\in{\Bbb P}(V)\times{\Bbb P}_5^0\times(\p2)^{n-2}\ \vert\\
&\text{$N_i:=v(S_i)$ is a node on $\{\ f=0\ \}$, $i=1,\dots,5$, $f(A_j)=0$, $j=3,\dots,n$}\ \}
\endalign
$$
with projection $\pi \colon X\to{\Bbb P}_5^0\times(\p2)^{n-2}$. Combining the proofs of Lemma~2.5.1 and Lemma 5.3.1 one obtains the following

\proclaim{Lemma 5.5.1}
There exists an open non--empty subset ${\Cal U}\subseteq(\p2)^{n-2}$ such that the scheme $X:=X'\cap({\Bbb P}(V)\times{\Bbb P}_5^0\times{\Cal U})$ is a subbundle of ${\Bbb P}(V)\times{\Bbb P}_5^0\times{\Cal U}$ with typical fibre $\p{12-n}$ over ${\Bbb P}_5^0\times{\Cal U}$.

Moreover the general point in $X$ corresponds to an irreducible sextic curve through the points $E_1,E_2$, containing exactly five nodes on 
$\widehat{C}$ and no other singularities.
\qed
\endproclaim

As above we conclude ${\Cal M}_{5,n}\approx X/T$.

\medbreak
Finally we are able to prove the following result.
\proclaim{Main Theorem for $g=5$} 
The moduli space ${\Cal M}_{5,n}$ is rational for $1\le n\le12$.
\endproclaim
\demo{Proof}
In all the cases $X$ is contained in a suitable $T$--linearized projective bundle over a rational base $B$ contained in the trivial $T$--equivariant projective bundle ${\Bbb P}(V)\times B$. Clearly there exists a $T$--invariant unisecant $L\subseteq {\Bbb P}(V)\times B$  (i.e. a divisor intersecting the general fibre in a hyperplane): the scheme $L\cap X$ is then a $T$--invariant unisecant on $X$, hence $X$ is $T$--equivariantly birationally equivalent to a $T$--equivariant vector bundle $\C^{R(n)}\times B$ where $R(1)=R(2)=10$ and $R(n)=12-n$ for $3\le n\le12$ as shown in Sections~5.3,~5.4,~5.5.

The action of $T$ on the base $B$ is in all cases almost free, thus ${\Cal M}_{5,n}\approx X/T\approx \C^{R(n)\times B}/T$ is a vector bundle over the base $B/T$ with $R(n)$--dimensional fibre due to the results of Section~4 of [Do]. It remains then to prove that the unirational variety $B/T$ is actually rational. To this purpose, we distinguish the cases $n=1,2$ and $n\ge3$.

If $n=1,2$ then $B$ is a projective linear representation of the connected 
solvable group $T$, hence $B/T$ is rational (see [Mi], [Vi]).

If instead $n\ge3$ then apply the results of [Do], Section~4, to 
the projective bundle 
${\Bbb P}_5^0\times(\p2)^{n-2}\to{\Bbb P}_5^0$. 
Once again, rationality follows and the proof is over.
\qed
\enddemo

\enddocument